\documentclass[preprint,12pt]{elsarticle}
\usepackage{latexsym, bm}

\usepackage{amssymb}
\usepackage{amsmath}

\journal{}
\def\th#1{\vspace{1mm}\noindent{\bf #1}\quad }
\begin{document}

\begin{frontmatter}



\title{The Goldstine-Weston theorem in random normed modules}


\author{Guang Shi\fnref{*}}
\fntext[*]{Supported by NNSF No. 10871016}

\address{LMIB and School of Mathematics and Systems Science, Beihang University, Beijing 100191, P.R. China }

\begin{abstract}
This article generalize the classical Goldstine-Weston theorem on normed spaces to one on random normed modules: the image of a random normed module $(E,\|\cdot\|)$ under the random natural embedding $J$ is dense in its double random conjugate space $E^{**}$ with respect to the $(\varepsilon,\lambda)$ weak star topology; and $J(E)$ is also dense in $E^{**}$  with respect to the locally $L^{0}$-convex weak star topology if $E$ has the countable concatenation property.

\end{abstract}


\end{frontmatter}



\section{Introduction}
\label{}

Random normed modules (briefly, $RN$ modules) is proved to be a proper and powerful tool in the study of conditional risk measures\cite{DMN,rece}. This successful application partly attributes to the systematic and deep development of the theory of $RN$ modules and random conjugate spaces\cite{exte,acha,bana,sepa,char,dual,rand,james}. As pointed out in \cite{bana}, the structure of an $RN$ module is an organic combination of the structure of a normed space with the structure of a measurable space. As a result, some classical theorem on normed spaces may still hold by the $RN$ modules, like James theory\cite{james}; and others may not hold universally for random normed modules unless they possess extremely simple stratification structure, like the Banach-Alaoglu theorem\cite{bana}. Exploring how the classical properties of normed spaces are preserved in $RN$ modules is important for further study and application of $RN$ modules.

The purpose of this paper is to establish similar results in $RN$ modules corresponding to the classical Goldstine-Weston theorem in normed spaces. Our results shows that the random weak star topologies on $RN$ modules process many similar properties to the classical weak star topology. Besides, one could also see the Goldstine-Weston theorem in $RN$ modules is exactly the same as the classical one under the $(\varepsilon,\lambda)$ topology; and slightly different from the classical one under the locally $L^{0}$-convex topology. This is due to that the $(\varepsilon,\lambda)$ topology is very natural and the locally $L^{0}$-convex topology is very strong.

The key step in our approach is that we generalize the classical Helly theorem characterizing the existence of solution for linear equations on a normed space to one characterizing the existence of solution for random linear equations on a random normed module. On one hand, we make full use of the recently developed theory of $RN$ modules in the prove of the Helly theorem in $RN$ modules; and on the other hand, the establishment of our final result depends much more on the Helly theorem in $RN$ modules than that in the normed spaces since the Banach-Alaoglu theorem is not universally available in $RN$ modules as mentioned above.

The paper proceeds as follows: section 2 gives some necessary notions and preliminaries. In section 3, we proved the helly theorem in random normed modules. And in section 4, we state and prove our main results.

\section{Preliminaries}
\label{}

Throughout this paper, $(\Omega,\mathcal{F},P)$ denotes a probability space, $N$ denotes the set of positive integers, and $K$ the scalar field $R$ of real numbers or $C$ of complex numbers.
Let $L^{0}(\mathcal{F},R)$ be the set of equivalence classes of real-valued
$\mathcal F$--measurable random variables on $\Omega$. Define the ordering $\leqslant$ on $L^{0}(\mathcal{F},\bar{R})$ by $\xi \leqslant\eta$ iff $\xi^{0}(\omega) \leqslant \eta^{0}(\omega)$ for $P$--almost all $\omega$ in $\Omega$
(briefly, a.s.), where $\xi^{0}$ and $\eta^{0}$ are arbitrarily chosen representatives of $\xi$ and $\eta$, respectively. It is well known from \cite{line} that $L^{0}(\mathcal{F},R)$ is a complete lattice in the sense that every subset with an upper bound has a supremum.
If $A$ is a subset with an upper bound of $L^{0}(\mathcal{F},R)$, we denote the supremum of $A$ by $\vee A$.\\

Besides, for any $A\in \mathcal{F}$ and $\xi ,\eta\in L^{0}(\mathcal{F},\bar{R})$, $``\xi >\eta$ on $A"$
means $\xi^{0}(\omega)>\eta^{0}(\omega)$ a.s. on $A$ for any chosen representative $\xi^{0}$ and $\eta^{0}$ of
$\xi$ and $\eta$, respectively. As usual, $\xi >\eta$ means $\xi \geqslant \eta$ and $\xi \neq \eta$; and $[\xi >\eta]$ denotes the equivalence class of the $\mathcal{F}$-measurable set $\{\omega\in \Omega~|~\xi^{0}(\omega)>\eta^{0}(\omega)\}$. Moreover, $L^{0}_{+}=\{\xi \in L^{0}(\mathcal{F},R)~|~ \xi \geqslant 0\}$ and $L^{0}_{++}=\{\xi \in L^{0}(\mathcal{F},R)~|~\xi > 0 ~ \mbox{on} ~ \Omega\}$.\\

And in the sequel of this paper we make the following convention: if $I_{A}$ denotes the characteristic function of an $\mathcal F$--measurable set $A$, then we use $\tilde{I}_{A}$ for its equivalence class in $L^{0}(\mathcal{F},K)$. Besides, for any $\xi \in L^{0}(\mathcal{F},K)$, $|\xi|$ and $\xi^{-1}$ respectively stand for the equivalence classes determined by the $\mathcal{F}$-measurable function $|\xi^{0}|:\Omega\rightarrow R$ defined  by $|\xi^{0}|(\omega)=|\xi^{0}(\omega)|$, $\forall \omega \in \Omega$ and $(\xi^{0})^{-1}$ defined by

$$
(\xi^{0})^{-1}(\omega)=\left \{
\begin{array}{ll}
(\xi^{0}(\omega))^{-1}, &\xi^{0}(\omega)\neq 0;\\
0, &\mbox{otherwise},
\end{array}
\right.
$$ where $\xi^{0}$ is an arbitrarily chosen representative of $\xi$. It is clear that $|\xi| \in L^{0}_{+}$
and $\xi \cdot \xi^{-1}=\tilde{I}_{\{\omega \in \Omega|\xi^{0}(\omega)\neq 0\}}$.\\

\th{Definition 2.1 \cite{rela}.}Let $E$ be a left module over the algebra $L^{0}(\mathcal{F},K)$. A countable concatenation of some
sequence $\{x_{n}~|~n\in N\}$ in $E$ with respect to some countable partition $\{A_{n}~|~n\in N\}$ of $\Omega$
is a formal sum $\Sigma_{n\in N}\tilde{I}_{A_{n}}x_{n}$. Moreover, a countable concatenation
$\Sigma_{n\in N}\tilde{I}_{A_{n}}x_{n}$ is well defined or $\Sigma_{n\in N}\tilde{I}_{A_{n}}x_{n}\in E$
if there is $x\in E$ such that $\tilde{I}_{A_{n}}x=\tilde{I}_{A_{n}}x_{n}$, $\forall n\in N$.
A subset $A$ of $E$ is called having the countable concatenation property if every countable
concatenation $\Sigma_{n\in N}\tilde{I}_{A_{n}}x_{n}$ with $x_{n}\in A$ for each $n\in N$
still belongs to $A$, namely $\Sigma_{n\in N}\tilde{I}_{A_{n}}x_{n}$ is well defined and there
exists $x\in A$ such that $x=\Sigma_{n\in N}\tilde{I}_{A_{n}}x_{n}$. And for a subset $M$ of $E$, $H_{cc}(M)$ denotes the countable catenation hull of $M$. \\

Suppose $M$ and $G$ are two nonempty subsets of an $L^{0}(\mathcal{F},K)$-module $E$. If $M$ and $G$
have the countable concatenation property and $M\cap G= \emptyset$, then $H(M,G)$ denotes the hereditarily disjoint stratification of $M$ and $G$\cite[Definition 3.14]{rela}.\\

\th{Definition 2.2 \cite{rela, base}.}An ordered pair $(E,\|\cdot\|)$ is called a random normed module over $K$ with base ($\Omega,\mathcal{F},P$) if $E$ is a left module over the algebra $L^{0}(\mathcal{F},K)$ and $\|\cdot \|$ is a mapping from $E$ to $L_{+}^{0}$ such that the following three axioms are satisfied:

\renewcommand{\labelenumi}{$($\arabic{enumi}$)$}
\begin{enumerate}

\item $\|x\|=0$ iff $x=\theta$(the null element of $E$);

\item $\|\xi x\|=|\xi|\|x\|$, $\forall \xi \in L^{0}(\mathcal{F},K)$ and $\forall x \in E$;

\item $\|x+y\|\leqslant \|x\|+\|y\|$, $\forall x,y \in E$.\\

\end{enumerate}

\th{Definition 2.3 \cite{rela, base}.}An ordered pair $(E,\langle\cdot,\cdot\rangle)$ is called a random inner product module (briefly, an $RIP$ module) over $K$ with base ($\Omega,\mathcal{F},P$) if $E$ is a left module over the algebra $L^{0}(\mathcal{F},K)$ and
$\langle\cdot,\cdot\rangle:E\times E\rightarrow L^{0}(\mathcal{F},K)$ satisfies the following statements:
\renewcommand{\labelenumi}{$($\arabic{enumi}$)$}
\begin{enumerate}
\item $\langle x,x\rangle\in L^{0}_{+}$ and $\langle x,x\rangle=0$ iff $x=\theta$;

\item $\langle x,y\rangle=\overline{\langle y,x\rangle}$, $\forall x,y \in E$ where $\overline{\langle y,x\rangle}$
denotes the complex conjugate of $\langle y,x\rangle$;

\item $\langle \xi x,y\rangle=\xi\langle x,y\rangle$, $\forall \xi \in L^{0}(\mathcal{F},K)$ and $\forall x,y \in E$;

\item $\langle x+y,z\rangle=\langle x,z\rangle+\langle y,z\rangle$,$\forall x,y,z\in E$.
\end{enumerate}
\noindent where $\langle x,y\rangle$ is called the random inner product between $x$ and $y$.\\

An $RIP$ module $(E,\langle\cdot,\cdot\rangle)$ is also an $RN$ module when $\|\cdot\|:E\rightarrow L^{0}_{+}$ is defined
by $\|x\|=\sqrt{\langle x,x\rangle}$, $\forall x\in E$.\\

\th{Example 2.4.}
Denote by $L^{0}(\mathcal{F},K^{n})$ the linear space of equivalence classes of $K^{n}$--valued $\mathcal F$--measurable functions on $\Omega$, where $n$ is a positive integer. Define
$\cdot:L^{0}(\mathcal{F},K)\times L^{0}(\mathcal{F},K^{n})\rightarrow L^{0}(\mathcal{F},K^{n})$ by
$\lambda \cdot x=(\lambda\xi_{1},\lambda\xi_{2},\cdots,\lambda\xi_{n})$
and
$\langle\cdot,\cdot\rangle:L^{0}(\mathcal{F},K^{n})\times L^{0}(\mathcal{F},K^{n})\rightarrow L^{0}(\mathcal{F},K)$
by $\langle x,y\rangle=\Sigma _{i=1}^{n}\xi_{i}\bar{\eta}_{i}$, \linebreak[4] for any $\lambda\in L^{0}(\mathcal{F},K)$ and $x=(\xi_{1},\xi_{2},\cdots,\xi_{n})$, $y=(\eta_{1},\eta_{2},\cdots,\eta_{n})\in L^{0}(\mathcal{F},K^{n})$.
It is easy to check that $(L^{0}(\mathcal{F},K^{n}),\langle\cdot,\cdot\rangle)$ is an $RIP$ module over $K$ with base ($\Omega,\mathcal{F},P$), and also an $RN$ module. Specially, $L^{0}(\mathcal{F},K)$ is an $RN$ module and $\|\lambda\|=|\lambda|$ for any $\lambda\in L^{0}(\mathcal{F},K)$.\\

\th{Definition 2.5 \cite{base}.}Let $(E,\|\cdot\|)$ be an $RN$ module over $K$ with base ($\Omega,\mathcal{F},P$). Then a linear operator $f$ from $E$ to $L^{0}(\mathcal{F},K)$ is called a $P$-a.e. bounded random linear functional on $E$ if there exists some $\xi$ in $L^{0}_{+}$ such that $|f(x)|\leqslant \xi \cdot \|x\|$, $\forall x \in S$.\\

Denote by $E^{*}$ the linear space of all $P$-a.e.~bounded random linear functionals on an $RN$ module $(E,\|\cdot\|)$ over $K$ with base ($\Omega,\mathcal{F},P$). Define
$\|\cdot\|^{\ast}: E^{*} \rightarrow L^{0}_{+}$ by
$\|f\|^{\ast}=\vee\{|f(y)|~|~y\in E$ and $\|y\|\leqslant 1,\}$ and
$\cdot:L^{0}(\mathcal{F},K)\times E^{\ast}\rightarrow E^{*}$ by $(\xi\cdot f)(x)=\xi \cdot (f(x))$,
$\forall \xi \in L^{0}(\mathcal{F},K)$, $\forall f\in E^{*}$ and $\forall x \in E$. Then it is easy to check that $(E^{*},\|\cdot\|^{*})$ is an $RN$ module over $K$ with base ($\Omega,\mathcal{F},P$). Such $(E^{*},\|\cdot\|^{*})$
is called the random conjugate space of $(E,\|\cdot\|)$\cite{base}.\\

Finally, this section is ended by introducing some useful topologies defined on $RN$ modules and their random conjugate spaces. A natural topology for an $RN$ modules $(E,\|\cdot\|)$ over $K$ with base ($\Omega,\mathcal{F},P$) is the $(\varepsilon,\lambda)$-topology\cite{base}, which is denoted by $\mathcal{T}_{\varepsilon,\lambda}$ in this paper. A subset $A$ of $E$ is $\mathcal{T}_{\varepsilon,\lambda}$ open iff for each $x\in E$, there exist two positive real numbers $\varepsilon,\lambda$ such that $\lambda<1$ and
$N_{x}(\varepsilon,\lambda)=\{y\in E~|~P([\|x-y\|<\varepsilon])>1-\lambda\}$ is included in $A$. And another stronger topology on $E$ is the locally $L^{0}$-convex topology\cite{DMN}, which is denoted by $\mathcal{T}_{c}$ in this paper. A subset $B$ is $\mathcal{T}_{c}$ open iff for each $x\in B$ there exists $\epsilon\in L^{0}_{++}$ such that $N_{x}(\varepsilon)=\{y\in E~|~\|y-x\|<\epsilon\}$ is included in $B$. Furthermore, similarly to the classical conjugate spaces, there are another two topologies for the random conjugate spaces $(E^{*},\|\cdot\|^{*})$ of $E$ besides the above two topologies, namely the $(\varepsilon,\lambda)$ weak star topology(denoted by $\sigma_{(\varepsilon,\lambda)}(E^{*},E)$) and the locally $L^{0}$-convex weak star topology(denoted by $\sigma_{c}(E^{*},E)$). If $\theta^{*}$ is the null element of $E^{*}$, then typical neighborhood systems of the null element $\theta^{*}$ of $E^{*}$ in $\sigma_{(\varepsilon,\lambda)}(E^{*},E)$ and $\sigma_{c}(E^{*},E)$ are the collection of the set
$N_{\theta^{*}}(x_{1},x_{2},\cdots,x_{n},\varepsilon,\lambda)=\{g\in E^{*}~|~P([|g(x_{i})|<\varepsilon])>1-\lambda,1\leqslant i\leqslant n\}$ and the collection of $N_{\theta^{*}}(x_{1},x_{2},\cdots,x_{n},\epsilon)=\{g\in E^{*}~|~|g(x_{i})|<\epsilon,1\leqslant i\leqslant n\}$ for all $n\in N$, $x_{1},x_{2},\cdots,x_{n}\in E$, $\varepsilon>0$, $0<\lambda<1$ and $\epsilon\in L^{0}_{++}$.
For details of these two topologies, and also terminologies such as $L^{0}$-convex, $L^{0}$-absorbent and $L^{0}$-balanced, we refer the readers to \cite{rela,DMN}.\\

\section{The Helly Theorem in random normed modules}
\label{}

\newtheorem{lem1}{Lemma}[section]
\begin{lem1}{
Let $(E,\|\cdot\|)$ be an $RN$ module over $R$ with base $(\Omega,\mathcal{F},P)$. Suppose $G$ and $M$ are two nonempty $L^{0}$-convex subsets of $E$ with the countable concatenation property and  the
$\mathcal{T}_{c}$ interior $G^{\circ}$ of $G$ is not empty. If $G\cap M=\emptyset$, then there exists $f \in E^{*}$ such that

$f(x)\leqslant f(y)$ on $H(G,M)$ for all $x \in G$ and $y\in M$

\noindent and

$f(x)<f(y)$ on $H(G,M)$ for all $x\in G^{\circ}$ and $y\in M$.

\noindent If $R$ is replaced by $C$, then the above statements still hold in the following way:

$(Ref)(x)\leqslant (Ref)(y)$ on $H(G,M)$ for all $x \in G$ and $y\in M$

\noindent and

$(Ref)(x)< (Ref)(y)$ on $H(G,M)$ for all $x \in G^{\circ}$ and $y\in M$.

\noindent Here $H(G,M)$ denotes the hereditarily disjoint stratification of $H$ and $M$,
and $(Ref)(x)=Re(f(x))$, $\forall x\in E$.
}\end{lem1}

\newproof{pf1}{Proof}
\begin{pf1}{
By \cite[Theorem 3.13]{rela} it follows that $P(H(G,M))>0$. First suppose that $H(G,M)=\Omega$.

Let $A=M-G=\{y-x~|~x\in G \;\mbox{and}\; y\in M\}$.
Clearly the $\mathcal{T}_{c}$ interior $A^{\circ}$ of $A$ is nonempty.
Define $B=\{z-x~|~x\in A\}$ for some fixed $z\in A^{\circ}$. It follows that $B$ is an $L^{0}$-convex subset of $E$;
and $B$ is also $L^{0}$-absorbent since $B$ contains a $\mathcal{T}_{c}$-neighborhood of $\theta$, where
$\theta$ is the null element of $E$. Moreover, it is easy to check that $H(\{z\},B)=H(G,M)=\Omega$.
Thus by \cite[Proposition 2.23 and Proposition 2.25]{DMN} the gauge function $p_{B}$ of $B$ is a random sublinear
function\cite{sepa} such that $p_{B}(z)\geqslant 1$ and $p_{B}(x)<1$ on $\Omega$ for any $x\in B^{\circ}$. Here $B^{\circ}$ is the $\mathcal{T}_{c}$ interior of $B$.

Define $g: \{\xi z~|~\xi\in L^{0}(\mathcal{F},R)\}\rightarrow L^{0}(\mathcal{F},R)$ by
$g(\xi z)=\xi p_{B}(z)$, $\forall \xi\in L^{0}(\mathcal{F},R)$. It is easy to verify that $g(\xi z)\leqslant p_{B}(\xi z)$,
$\forall \xi\in L^{0}(\mathcal{F},R)$. Thus by the Hahn-Banach extension theorem in $RN$ module\cite[Theorem2.8]{rela},
there exists a random linear functional $f: E\rightarrow L^{0}(\mathcal{F},R)$ such that $f$ extends $g$
and $f(x)\leqslant p_{B}(x)$, $\forall x\in E$.

Notice that $f(x)\leqslant p_{B}(x)<1$ on $\Omega$ and $f(x)=-f(-x)\geqslant -p_{B}(x)>-1$ on $\Omega$ for any
$x\in B^{\circ}$. It follows that
$$f(B^{\circ})\subset \{\xi\in L^{0}(\mathcal{F},R)|\; |\xi|<1\; \mbox{on}\; \Omega\}.$$
Hence $f$ is a $\mathcal{T}_{c}$ continuous module homomorphism, i.e. $f\in E^{*}$. Moreover, since $z-(y-x)\in B$ for
any $x\in G$ and $y\in M$, thus $$
f(x-y)=f(z-(y-x))-f(z)\leqslant 1-p_{B}(z) \leqslant 0,$$
i.e. $f(x)\leqslant f(y)$. Likewise, since $z-(y-x)\in B^{\circ}$ for
any $x\in G^{\circ}$ and $y\in M$, thus $$
f(x-y)=f(z-(y-x))-f(z)<1-p_{B}(z) \leqslant 0 \;\mbox{on}\;\Omega,$$
i.e. $f(x)<f(y)$ on $\Omega$. Hence $f$ is the required random functional.

If $H(G,M)\neq\Omega$, let $\Omega^{\prime}=H(G,M)$,
$\mathcal{F}^{\prime}=\Omega^{\prime}\cap \mathcal{F}=\{\Omega^{\prime}\cap F~|~F\in \mathcal{F}\}$
and $P^{\prime}:\mathcal{F}^{\prime}\rightarrow [0,1]$ be defined by
$P^{\prime}(\Omega^{\prime}\cap F)=P(\Omega^{\prime}\cap F)/P(\Omega^{\prime})$.
Take $E^{\prime}=\tilde{I}_{\Omega^{\prime}}E$, $M^{\prime}=\tilde{I}_{\Omega^{\prime}}M$,
$G^{\prime}=\tilde{I}_{\Omega^{\prime}}G$ and consider $(E^{\prime},\|\cdot\|_{E^{\prime}})$ as
an $RN$ module with base $(\Omega^{\prime},\mathcal{F}^{\prime},P^{\prime})$. Then $M^{\prime}$
and $G^{\prime}$ satisfy the above condition, so there exists an $f^{\prime}\in (E^{\prime})^{*}$
such that

$f^{\prime}(x)\leqslant f^{\prime}(y)$ on $\Omega^{\prime}$ for all $x \in G^{\prime}$ and $y\in M^{\prime}$

\noindent and

$f^{\prime}(x)< f^{\prime}(y)$ on $\Omega^{\prime}$ for all $x \in (G^{\prime})^{\circ}$ and $y\in M^{\prime}$.

\noindent
By \cite[Theorem 2.10]{rela} $f^{\prime}$ has an extension $f\in E^{*}$, which meets our requirement.

Finally, if $R$ is replaced by $C$, the result follows immediately by noticing that every
$RN$ module over $C$ is also an $RN$ module over $R$ and

$f(x)=(Ref)(x)-i(Ref)(ix)$,
$\forall f \in E^{*}$ and $\forall x\in E$.
$\Box$
}\end{pf1}

\newtheorem{thm1}[lem1]{Theorem}
\begin{thm1}{
Suppose $(E,\|\cdot\|)$ is an $RN$ module over $C$ with base $(\Omega,\mathcal{F},P)$ and $E$ has
the countable concatenation property. $f_{1},f_{2},\cdots,f_{n} \in E^{*}$,
$\xi_{1}, \xi_{2},\cdots ,\xi_{n} \in L^{0}(\mathcal{F},C)$ and
$\beta \in L_{+}^{0}$.
For any $\varepsilon \in L^{0}_{++}$, there exists $x_{\varepsilon} \in E$ which satisfies:\\

(1) $f_{i}(x_{\varepsilon})=\xi _{i}$, $i=1,2,\cdots ,n$;

(2) $\|x_{\varepsilon}\|\leqslant \beta+\varepsilon$\\

\noindent iff

$$|\sum_{k=1}^{n}\lambda_{k}\xi_{k}|\leqslant \beta\|\sum_{k=1}^{n} \lambda_{k}f_{k}\|$$

\noindent holds for arbitrary $\lambda_{1},\lambda_{2},\cdots \lambda_{n}\in L^{0}(\mathcal{F},C)$.
}\end{thm1}

\newproof{pf2}{Proof}
\begin{pf2}{

Necessity is obvious, it remains to prove sufficiency.

Let $S=\{\Sigma_{i=1}^{n}\zeta_{i}f_{i}~|~\zeta_{i}\in L^{0}(\mathcal{F},C), 1\leqslant i\leqslant n\}$, then $S$ is a finitely generated $L^{0}(\mathcal{F},C)$ modules. By \cite[Theorem 3.1]{fini}, there exists a partition $\{A_{0},A_{1},\cdots,A_{n}\}$ of $\Omega$ to ${\cal F}$ such that $\tilde{I}_{A_{i}}S$ is a quasi-free stratification of rank $i$ of $S$ for each $i$ which satisfies $0\leqslant i\leqslant n$ and $P(A_{i})>0$. Let $\{g_{j}\in \tilde{I}_{A_{i}}S~|~1\leqslant j\leqslant i\}$ be a basis for $\tilde{I}_{A_{i}}S$ for some $i$ such that $1\leqslant i\leqslant n$ and $P(A_{i})>0$. Suppose $g_{j}=\Sigma_{k=1}^{n}\zeta_{kj}f_{k}$, $1\leqslant j\leqslant i$ and $\zeta_{kj}\in \tilde{I}_{A_{i}}L^{0}(\mathcal{F},C)$. Let $\gamma_{j}=\Sigma_{k=1}^{n}\zeta_{kj}\xi_{k}$, then
$$|\sum_{j=1}^{i}\lambda_{j}\gamma_{j}|=|\sum_{j=1}^{i}\sum_{k=1}^{n}\lambda_{j}\zeta_{kj}\xi_{k}|\leqslant \beta\|\sum_{j=1}^{i}\sum_{k=1}^{n}\lambda_{j}\zeta_{kj}f_{k}\|=\beta\|\sum_{j=1}^{i}\lambda_{j}g_{j}\|.$$
If there exists $x_{A_{i}}\in \tilde{I}_{A_{i}}E$ such that $g_{j}(x_{A_{i}})=\gamma_{j}$ for $1\leqslant j\leqslant i$, then $f_{k}(x_{A_{i}})=\tilde{I}_{A_{i}}\xi_{k}$ for $1\leqslant k\leqslant n$. Actually, suppose $\tilde{I}_{A_{i}}f_{k}=\sum_{j=1}^{i}\eta_{jk}g_{j}$ for $1\leqslant k\leqslant n$ and $\eta_{jk}\in \tilde{I}_{A_{i}}L^{0}(\mathcal{F},C)$, then
$$|\tilde{I}_{A_{i}}\xi_{k}-\sum_{j=1}^{i}\eta_{jk}\gamma_{j}|\leqslant\beta\|\tilde{I}_{A_{i}}f_{k}-\sum_{j=1}^{i}
\eta_{jk}g_{j}\|=0,$$
i.e. $\tilde{I}_{A_{i}}\xi_{k}=\sum_{j=1}^{i}\eta_{jk}\gamma_{j}$. Hence
$$f_{k}(x_{A_{i}})=\sum_{j=1}^{i}\eta_{jk}g_{j}(x_{A_{i}})=\sum_{j=1}^{i}\eta_{jk}\gamma_{j}=\tilde{I}_{A_{i}}\xi_{k}.$$

Thus no lose of generality, suppose
$\{f_{1},f_{2},\cdots,f_{n}\}$ is $L^{0}(\mathcal{F},C)$-independent.
Define $T:E\rightarrow L^{0}(\mathcal{F},C^{n})$
by $Tx=(f_{1}(x),f_{2}(x),\cdots ,f_{n}(x))$, $\forall x \in E$.
Obviously $T(E)$ is a submodule of $L^{0}(\mathcal{F},C^{n})$ with the countable concatenation property.
If $T(E)\neq L^{0}(\mathcal{F},C^{n})$, by \cite[Corollary 4.3]{fini} there exists
$z=(\eta_{1},\eta_{2},\cdots,\eta_{n})\in L^{0}(\mathcal{F},C^{n})$ such that
$$
(\sum_{k=1}^{n} \bar{\eta}_{k}f_{k})(x)=\sum_{k=1}^{n} \bar{\eta}_{k}f_{k}(x)=\langle T(x),z\rangle=0, \forall  x\in E.
$$
This contradicts with the $L^{0}(\mathcal{F},C)$-independence of $\{f_{1},f_{2},\cdots,f_{n}\}$. Hence $T(E)=L^{0}(\mathcal{F},C^{n})$.

Suppose $x_{1},x_{2},\cdots,x_{n}\in E$ such that $Tx_{i}=(\eta_{1}^{i},\eta_{2}^{i},\cdots,\eta_{n}^{i})$,
$\eta_{i}^{i}=1$ and $\eta_{j}^{i}=0(i\neq j)$ for $1\leqslant j \leqslant n$, $1\leqslant i \leqslant n$.
Let $\gamma=\vee_{i=1}^{n}\|x_{i}\|$, then clearly $\gamma>0$ on $\Omega$. If
$y=(\alpha_{1},\alpha_{2},\cdots,\alpha_{n}) \in L^{0}(\mathcal{F},C^{n})$ and
$\|y\|\leqslant(\beta+\varepsilon)n^{-1}\gamma^{-1}$ for some fixed $\varepsilon\in L^{0}_{++}$,
then $T(\sum_{i=1}^{n} \alpha_{i}x_{i})=y$ and $$
\|\sum_{i=1}^{n} \alpha_{i}x_{i}\| \leqslant \sum_{i=1}^{n} |\alpha_{i}|\|x_{i}\|\leqslant(\beta+\varepsilon)n^{-1}\gamma^{-1}\sum_{i=1}^{n}\|x_{i}\|\leqslant\beta+\varepsilon.$$
Let $\bar{B}_{\beta+\varepsilon}=\{x\in E|\;\|x\|\leqslant\beta+\varepsilon\}$, the above argument shows that
$T(\bar{B}_{\beta+\varepsilon})$ contains a $\mathcal{T}_{c}$-open neighborhood
$\{y\in L^{0}(\mathcal{F},C^{n})~|~ \|y\|<(\beta+\varepsilon)n^{-1}\gamma^{-1}\}$ of the null
element of $L^{0}(\mathcal{F},C^{n})$. It is easy to see that $T(\bar{B}_{\beta+\varepsilon})$
is also an $L^{0}$-convex subset with the countable concatenation property.

If the hypothesis does not hold, i.e. there exists $\varepsilon\in L^{0}_{++}$ such that
$p\triangleq(\xi_{1}, \xi_{2},\cdots ,\xi_{n})\notin T(\bar{B}_{\beta+\varepsilon})$.
Thus by Lemma 3.1 there exists $f\in L^{0}(\mathcal{F},C^{n})^{*}$ such that
$(Ref)(y)\leqslant (Ref)(p)$ on $H(\{p\},T(\bar{B}_{\beta+\varepsilon}))$ for all
$y \in T(\bar{B}_{\beta+\varepsilon})$ and
$(Ref)(p)>(Ref)(0)=0$ on $H(\{p\},T(\bar{B}_{\beta+\varepsilon}))$.
For any fixed $y\in T(\bar{B}_{\beta+\varepsilon})$, let $\xi=|f(y)|(f(y))^{-1}$, then
clearly $\xi y\in T(\bar{B}_{\beta+\varepsilon})$. Moreover,
$$|f(y)|=f(\xi y)=(Ref)(\xi y)\leqslant (Ref)(p)\leqslant |f(p)|$$

\noindent on $H(\{p\},T(\bar{B}_{\beta+\varepsilon}))$.

By Riesz's representation theorem in $RIP$ module
\cite[Theorem 4.3]{rela}, there exists $y_{0}=(\lambda_{1},\lambda_{2},\cdots \lambda_{n})\in L^{0}(\mathcal{F},C^{n})$
such that $f(y)=\langle y,y_{0}\rangle$, $\forall y\in L^{0}(\mathcal{F},C^{n})$. Hence
$$
|\sum_{k=1}^{n} \bar{\lambda}_{k}f_{k}(x)|=|f(Tx)|\leqslant |f(p)|=|\sum_{k=1}^{n} \bar{\lambda}_{k}\xi_{k}|
$$
\noindent on $H(\{p\},T(\bar{B}_{\beta+\varepsilon}))$, $\forall x\in \bar{B}_{\beta+\varepsilon}$.
Thus
$$
(\beta+\varepsilon)\|\sum_{k=1}^{n} \bar{\lambda}_{k}f_{k}\|=\bigvee_{x\in \bar{B}_{\beta+\varepsilon}}
|\sum_{k=1}^{n} \bar{\lambda}_{k}f_{k}(x)|\leqslant |\sum_{k=1}^{n} \bar{\lambda}_{k}\xi_{k}|.$$

\noindent on $H(\{p\},T(\bar{B}_{\beta+\varepsilon}))$.

Since $|f(p)|=|\langle p,y_{0}\rangle |>0$ on $H(\{p\},T(\bar{B}_{\beta+\varepsilon}))$, it follows
$\|y_{0}\|>0$ on $H(\{p\},T(\bar{B}_{\beta+\varepsilon}))$. Thus $\|\sum_{k=1}^{n} \bar{\lambda}_{k}f_{k}\|\neq 0$
on $H(\{p\},T(\bar{B}_{\beta+\varepsilon}))$ by the $L^{0}(\mathcal{F},C)$-independence of $\{f_{1},f_{2},\cdots,f_{n}\}$. Hence
$$
\beta\|\sum_{k=1}^{n} \bar{\lambda}_{k}f_{k}\|<|\sum_{k=1}^{n} \bar{\lambda}_{k}\xi_{k}|
$$
\noindent on $H(\{p\},T(\bar{B}_{\beta+\varepsilon}))$, which contradicts with the assumption.$\Box$
}\end{pf2}

\newdefinition{rmk}[lem1]{Remark}
\begin{rmk}{It is necessary to require $E$ to have the countable concatenation property in Theorem 4.2, otherwise
the result may not hold. Here is an example. Let $\Omega=[0,1]$, ${\cal F}$ be the collection of all Lebesgue measurable subsets of $[0,1]$ and $P$ the Lebesgue measure on $[0,1]$. Suppose $M=\{\tilde{I}_{[2^{-(n+1)},2^{-n}]}~|~n\in N\}$ and $E=\{\sum_{i=1}^{n}\xi_{i}x_{i}~|~\xi_{i}\in L^{0}(\mathcal{F},C),x_{i}\in M,1\leqslant i\leqslant n$ and $n\in N\}$. Clearly $E$ is a submodule of
$L^{0}(\mathcal{F},C)$ without the countable concatenation property. Define $\|\cdot\|:E\rightarrow L^{0}_{+}$ by $\|\eta\|=|\eta|$, $\forall\eta\in E$, then
$(E,\|\cdot\|)$ is also an $RN$ module over $C$ with base $(\Omega,\mathcal{F},P)$. If $f \in E^{*}$ is defined
by $f(\eta)=\eta$, $\forall\eta \in E$ and take $\xi=\tilde{I}_{\Omega}$, $\beta=\tilde{I}_{\Omega}$. Then clearly $|\lambda\xi|\leqslant \beta\|\lambda f\|$(in fact, $|\lambda\xi|= \beta\|\lambda f\|$), $\forall\lambda \in L^{0}(\mathcal{F},C)$. But there does not exist
any $x\in E$ such that $f(x)=\xi$.

}\end{rmk}

\section{The Goldstine-Weston theorem in $RM$ modules}
\label{}

\newtheorem{lem2}{Lemma}[section]
\begin{lem2}{Let $(E,\|\cdot\|)$ be an $RN$ module over $C$ with base $(\Omega,\mathcal{F},P)$. If $E$ has the countable concatenation property, then the unit ball $E^{*}(1)=\{f\in E^{*}~|~\|f\|\leqslant1\}$ of the random conjugate space $E^{*}$ of $E$ is closed with respect to both $\sigma_{(\varepsilon,\lambda)}(E^{*},E)$ and $\sigma_{c}(E^{*},E)$.

}\end{lem2}

\newproof{pf3}{Proof}
\begin{pf3}{Since $\sigma_{c}(E^{*},E)$ is stronger than $\sigma_{(\varepsilon,\lambda)}(E^{*},E)$, it only needs to prove $E^{*}(1)$ is $\sigma_{(\varepsilon,\lambda)}(E^{*},E)$ closed. Suppose $g\in E^{*}$ and $g\notin E^{*}(1)$, then there exist $A\in \mathcal{F}$ and a positive real number $\delta$ such that $P(A)>0$ and $\|g\|>1+\delta$ on $A$. By Theorem 3.2 there is an $x\in E$ such that $\|x\|\leqslant(1+\delta/2)\tilde{I}_{A}$ and $g(x)=\tilde{I}_{A}\|g\|$. Then the $\sigma_{(\varepsilon,\lambda)}(E^{*},E)$ neighborhood $N_{g}(x,\delta/2,P(A)/2)$ of $g$ is disjoint with $E^{*}(1)$, which completes the proof.$\Box$

}\end{pf3}

\newtheorem{thm2}[lem2]{Theorem}
\begin{thm2}{Let $(E,\|\cdot\|)$ be an $RN$ module over $C$ with base $(\Omega,\mathcal{F},P)$ with the countable concatenation property, $J$ denote the random natural embedding $E\rightarrow E^{**}$ which is defined by $J(x)(g)=g(x)$ for any $g\in E^{*}$ and any $x\in E$. Then the $\sigma_{c}(E^{*},E)$ closure of $J(E)$ is $E^{**}$.
}\end{thm2}

\newproof{pf4}{Proof}
\begin{pf4}{The result will be established immediately once we proved that the $\sigma_{c}(E^{*},E)$ closure of $J(E(1))$ is $E^{**}(1)$. Suppose $\mathfrak{l}\in E^{**}(1)$, $f_{1},f_{2},\cdots,f_{n}\in E^{*}$ and $\varepsilon\in L^{0}_{++}$. Let $\gamma=\bigvee_{i=1}^{n}\|f_{i}\|\vee \tilde{I}_{\Omega}$. By Theorem 3.2, there exists $x_{0}\in E$ such that $\|x_{0}\|\leqslant1+2\varepsilon\|\gamma\|^{-1}$ and $f_{i}(x_{0})=\mathfrak{l}(f_{i})$ for $1\leqslant i\leqslant n$. Let $x=\gamma(\gamma+\varepsilon/2)^{-1}x_{0}$, then $\|x\|\leqslant1$ and
$$|(J(x)(f_{i})-\mathfrak{l}(f_{i}))|=|f_{i}(x-x_{0})|\leqslant\varepsilon\|x\|/2<\varepsilon$$
for $1\leqslant i\leqslant n$. Thus $J(x)$ belongs to the $\sigma_{c}(E^{*},E)$ neighborhood $\{\mathfrak{h}\in E^{**}~|~|(\mathfrak{h}-\mathfrak{l})(f_{i})|<\varepsilon,~i=1,2,\cdots,n\}$. Hence $J(E(1))$ is $\sigma_{c}(E^{*},E)$ dense in $E^{**}(1)$, which completes the proof.$\Box$

}\end{pf4}

\newtheorem{thm3}[lem2]{Theorem}
\begin{thm3}{Suppose $(E,\|\cdot\|)$ is an $RN$ module over $C$ with base $(\Omega,\mathcal{F},P)$, then the $\sigma_{(\varepsilon,\lambda)}(E^{*},E)$ closure of $J(E)$ is $E^{**}$.
}\end{thm3}

\newproof{pf5}{Proof}
\begin{pf5}{Suppose $E_{cc}=H_{cc}(E)$ and define $\|\cdot\|_{cc}:E_{cc}\rightarrow L_{+}^{0}$ by
$\|x\|_{cc}=\sum_{n\in N}\tilde{I}_{A_{n}}\|x\|_{n}$ for any $x=\sum_{n\in N}\tilde{I}_{A_{n}}x_{n}$ in $E_{cc}$, where $\{A_{n}~|~n\in N\}$ is a countable partition of $\Omega$ to $\mathcal{F}$ and $x_{n}\in E$ for $n\in N$. It is easy to check that $E_{cc}^{**}=E^{**}$. By the above theorem, $J(E_{cc})$ is dense in $E^{**}$ with respect to $\sigma_{(\varepsilon,\lambda)}(E^{*},E)$ since $\sigma_{c}(E^{*},E)$ is stronger than $\sigma_{(\varepsilon,\lambda)}(E^{*},E)$. Combine the fact that $J(E)$ is $\mathcal{T}_{\varepsilon,\lambda}$ dense in $J(E_{cc})$ with that $\mathcal{T}_{\varepsilon,\lambda}$ is stronger than  $\sigma_{(\varepsilon,\lambda)}(E^{*},E)$, it follows our desired result.$\Box$

}\end{pf5}











\end{document}